\documentclass[11pt]{myj}

\usepackage{amsmath}
\usepackage{amsfonts, amssymb, amsthm}
\newcommand{\tauG}{\tau(G)}
\newcommand{\Id}{\mathrm{Id}}

\newcommand{\Lap}{\Delta}

\newcommand{\cG}{{\mathcal{G}}}

\newcommand{\cT}{{\mathcal{T}}}

\newcommand{\reals}{\mathbb{R}}
\newcommand{\sgn}{\mathop{\mathrm{sgn}}}

\begin{document}

\title{Extremal metrics on graphs I}
\author{Dmitry Jakobson \and Igor Rivin}

\address{First Author: Mathematics Department, University of
Chicago.\\
 Second author: Mathematics Department, University of
Manchester and
  Mathematics Department, Temple University}

\email{rivin@ma.man.ac.uk\\jakobson@math.uchicago.edu}

\thanks{I.~Rivin would like to thank the \'Ecole Polytechnique for
its hospitality during the preparation of this paper. D.~Jakobson was 
partially supported by the NSF}

\date{\today}

\keywords{graphs, extremal graph theory, deformation theory, uniformization}

\begin{abstract}
We define a number of natural (from geometric and combinatorial
points of view) deformation spaces of valuations on finite graphs,
and study functions over these deformation spaces. These functions
include both direct metric invariants (girth, diameter), and spectral 
invariants (the determinant of the Laplace
operator, or complexity; bottom non-zero eigenvalue of the Laplace
operator). We show that almost all of these functions are,
surprisingly, convex, and we characterize the valuations
extremizing these invariants.
\end{abstract}

\maketitle

\newtheorem{prop}{Proposition}
\newtheorem{defn}{Definition}
\newtheorem{lemma}{Lemma}
\newtheorem{thm}{Theorem}
\newtheorem{fact}{Fact}
\newtheorem{cor}{Corollary}
\newtheorem{remark}{Remark}
\newtheorem{ex}{Example}
\newtheorem{question}{Question}
\newtheorem{obs}{Observation}
\newcommand{\tr}{\mathrm{tr\:}}
\newcommand{\spec}{\mathrm{spec\:}}

\section*{Introduction}
There is a vast literature on the subject of \textit{extremal graph
  theory}. There, the general approach is to consider a natural
  invariant (invariant with respect to isomorphism) of graphs, and to
  try to understand which graphs make the invariant as big as
  possible, subject to (presumably natural) constraints. Examples of
  such invariants are:

Girth -- the length of the shortest cycle;

Diameter -- the greatest distance between a pair of vertices;

Tree number -- the number of spanning trees,

and some closely related spectral invariants: the ``determinant of
the Laplacian'', the smallest positive eigenvalue of the
Laplacian, and so on.

It is expected that graphs which are ``good'' with respect to any one
of these invariant will be good with respect to the others (where by
``good'', we mean that the graph is either extremal, or close to it),
and will have other (\textit{a priori} unsuspected) nice combinatorial
properties.

Extremal graph theory is a rather difficult subject, largely due to
its intrinsically combinatorial nature 
(arguably it is this difficulty which attracts most of the
practitioners).

A seemingly not very closely related subject is that of
differential geometry. One of its central areas is that of
``uniformization'', or ``optimal geometry''. There, we are often
given a \textit{fixed} topological space, and we try to find a
metric on this space which maximizes some invariant. The actual
invariants studied are very often similar to those mentioned above
for graphs. The \textit{motivation}, on the other hand, is
sometimes the same as that of extremal graph theory, but sometimes
there is an additional factor: it is hoped that the extremal
metrics would give a \textit{canonical} representation of the
topological space, which renders its topological properties more
transparent (for example, the study of the topology of the sphere
would be much more difficult if we did not have its standard
``round'' representation at our disposal).

Our motivation stems from both the areas sketched above: we would
like to get canonical representations of graphs, but we have other
concerns as well. First of all, the space of edge valuations of a
given finite graph is a much simpler space than the space of
metrics on a given topological space. Thus, we hope that the
answers to our questions will be technically simpler than the
corresponding differential-geometric results, but that the model
is sufficiently rich to suggest what one might expect. By the same
token, the space of edge valuations on a fixed graph is a much
simpler space than the (discrete) space of graphs, though the
latter is naturally embedded in the former. We thus hope to get
insight into problems in extremal graph theory as well.

\subsection{Outline of the paper} We set up the basic deformation
spaces and announce the main convexity results in Section
\ref{found}. We set up the girth problem in Section \ref{girth},
and characterize the extremal valuations in Section \ref{girth2}. 

We define the basic matrices and operators we are working with,
and show the convexity of the bottom eigenvalue and the complexity
in Section \ref{spectr1}. We characterize valuations extremal for
complexity (or ``determinant of laplacian'' in Section
\ref{maxtreesec}, and valuations extremal for the bottom
eigenvalue in Section \ref{maxevsec}. Finally, in Section
\ref{constgrad} we analyse completely those graphs which are
extremal for $\lambda_1$, under the additional assumption (which
turns out to be very strong), that $\lambda_1$ appears without
multiplicity.

\section{The foundations}\label{found}

We will always consider a fixed \textit{finite simple graph} $G$.
We will consider the following deformation spaces of edge
valuations on $G$:

\begin{equation*}
P(G) = \left\{f: E(G) \rightarrow \mathbb{R}^+ \quad | \quad
\sum_{e
    \in E(G)} f(e) = |E(G)|\right\},
\end{equation*} -- the space of all edge valuations of
    $G$.

\begin{equation*}
T(G) = \left\{f\in P(G) \quad | \quad \sum_{\text{$e_i$ incident to
      $v$}}  = d_v, \forall v \in V(G)\right\}
\end{equation*}

\begin{equation*}
C(G) = \left\{f\in P(G) \ | \ \exists g:V(G) \rightarrow
  \mathbb{R}^+,\ \text{such that}\quad f(e_{vw}) = g(v) + g(w)\right\}
\end{equation*}

The letters $T$ and $C$ are meant to suggest Teichm{\"u}ller space and
conformal deformation space respectively.

All three spaces have a natural linear structure, which we will use
without further comment.

We will look at the variation of following invariants (defined
below) over the above-described deformation spaces: girth $g$,
bottom positive eigenvalue $\lambda_1$ of the
Laplacian $\Delta$ and $\log \det^* \Delta$.

The first striking observation about these invariants is the
following:

\begin{thm}
The quantities $-g$, $-\lambda_1$, $-\log \det^* \Delta$ are
\emph{convex} on $P(G)$ (and hence on its linear subspaces $C(G)$ and
$T(G)$).
\end{thm}

(The proofs of these results are spread out through this paper:
The convexity of girth is given in Section 
\ref{girth}; 
the convexity of $-\log \det^* \Delta$ -- by
Theorem \ref{lggconc}, and the convexity of $\lambda_1$ is
outlined in Section \ref{varspec}.)

\begin{remark}
It can also be shown that the ``topological entropy of the
geodesic flow'', defined in terms of a different deformation of
the adjacency matrix, is also convex. This is done in the article
\cite{Riv99} by the second author. \end{remark}
 The convexity has
far-reaching consequences. To wit, for every invariant $I \in \{g,
 \lambda_1, -\log \det^* \Delta\}$, and for each deformation
space $D \in \{P, C, T\}$ there is a \textit{unique}
\textit{canonical} edge valuation $G_D^I$ maximizing the
invariant. A natural question is one of the characterization of
these critical valuations, and of understanding the relationship
between the various $G_D^I$ for the different choices of $I$ and
$D$.

Some properties follow immediately from the convexity, in
particular:
\begin{obs}
If $G$ possesses a group $\Gamma(G)$ of automorphisms, then the
weights of the critical points $G_D^I$ are invariant by these
symmetries, thus, if the automorphism group of $G$ is edge transitive,
then $G_D^I$ are all equal (independently of invariant and deformation
space), and are given by the constant weighing on the edges. If the
automorphism group of $G$ acts vertex-transitively, then $G_C^I$ is
given by the constant weighing.
\end{obs}

\begin{remark}
A large class of graphs the automorphism group of which is vertex-
but not edge- transitive is given by the Cayley graphs of finite
groups.
\end{remark}

For graphs not known \textit{a priori} to be symmetric, the
supposition that the unweighted graph $G$ is critical for one of
the invariants, implies strong symmetry properties. For example,
if $G$ is maximal for $\log \det^* \Delta$, then there is the same
number of spanning trees through every edge of $G$ ($G$ is
\textit{equiarboreal} in the terminology of Godsil). If $G$ is the
maximum for $\lambda_1$ then, with rare exceptions, $\lambda_1$
occurs with multiplicity in the spectrum of $\Delta(G)$. If $G$ is
maximal for girth, then every edge of $G$ is contained in a
shortest cycle (the precise somewhat stronger statement is the
content of Theorem \ref{omnibus}).

\section{Girth}
\label{girth}
 The ``direct'' (girth) and
``spectral'' invariants are somewhat different conceptually. First
we remind the reader that the \textit{length} of a path in a
weighted graph $G$ is the sum of the weights of the edges in the
path.  The \textit{girth} $\gamma(G)$ is the length of the
shortest cycle in $G$. The \emph{distance} between two vertices of
$G$ is the length of the shortest path connecting them; the {\em
diameter} $D(G)$ is equal to the largest such distance.  Thus, if
$\mathfrak{C}$ is the set of all cycles of $G$, then the girth is
given by:
\begin{equation*}
g=\min_{C\in \mathfrak{C}}\sum_{e\in C} f(e),
\end{equation*}
where $f(e)$ is the valuation of the edge $e$. Note that each of the
terms $\sum_{e\in C} v(e)$ is a linear function of the valuation $f$,
and hence we have the immediate
\begin{thm}
The girth $g(f)$ is a concave function on $P(G)$.
\end{thm}

\begin{proof}
This follows from the observation that the minimum of a collection
of concave (in particular linear) functions is concave. We leave
the proof as an exercise to the interested reader.
\end{proof}

For two vertices $u,v$ of $G$ denote by
$\Lambda(u, v)$  the set of all paths $\lambda$ in $G$
connecting $u$ and $v$ (we can assume without loss of generality that
the paths are not self-intersecting, to make sure that the number of 
paths considered is finite). Thus, 
\begin{equation*} 
d(u, v) = \min_{\lambda \in \Lambda} \sum_{e\in \lambda} f(e). 
\end{equation*} 
The diameter of $G$ is thus given by:
\begin{equation*}
D(G) = \max_{u, v \in V(G)} d(u, v).
\end{equation*}
Note that the diameter is \emph{not} \textit{a priori} convex, due
to the additional maximum, though some of the methods we use for
girth can be brought to bear on the diameter question as well.

\section{Spectral invariants}
\label{spectr1} Let $G$ be a graph with $n$ vertices and $m$ edges
(we denote the set of such graphs by $\cG_{m,n}$), and $f: E(G)
\rightarrow \mathbb{R}$ be a valuation on the edges of $G$ (in the
sequel the valuations are always assumed positive, but this is not
essential for the definitions below). The {\em adjacency} matrix
$A(G)$ (or just $A$, when no ambiguity is possible) of a graph
$G\in \cG_{m,n}$ is a square matrix of size $n$ where $A_{i,j}$ is
weight of the edge joining the vertices $v_i$ and $v_j$ if there
such an edge, and $0$ otherwise. We always consider loopless
graphs, so $A_{ii}=0$.

The \textit{nearest neighbor  Laplacian} $\Lap$ acts on functions
on the set $V(G)$ of the vertices of $G$: given $g:V(G)\to\reals$,
$\Lap (g)(v)\ =\ \sum_{w\sim v} f(vw)(g(v)-g(w))$. (We
exceptionally use $vw$ to denote the edge joining $v$ to $w$).
 Let $\Lap =\{L_{ij}\}$ be the
matrix of $\Lap(G)$ of a (not necessarily simple) graph $G$; then $L_{ii}$ is
the degree of the vertex $v_i$, and $-L_{ij}$ is the number of the edges
joining $v_i$ and $v_j\neq v_i$ (equal to $0$ or $1$ for simple graphs).
For $k$-regular graphs $\Lap=k\Id -A$.

Let $A$ be the adjacency matrix of a weighted graph $G$ with $n$
vertices, let $\delta$ be the maximal degree of a vertex in $G$, and
let its spectrum (in the decreasing order) be given by
\begin{equation}\label{adj:spectrum}
   \delta\geq\mu_1 >\mu_2\geq\mu_3\geq\ldots\geq\mu_n\geq (-\delta). 
\end{equation}
The spectrum of $\Lap(G)$ is 
$0=\lambda_0<\lambda_1\leq\ldots\leq\lambda_{n-1}$.  
For $k$-regular graphs, $\delta=k=\mu_1$ and  
$\lambda_j=k-\mu_{j+1}, j=0,1,\ldots,n-1$.   Note that the
$0$ is always in the spectrum of $\Lap(G)$ independently of the
weighing on the edges, and, furthermore, as long as the weighing
is strictly positive, and the graph $G$ is connected, the
eigenspace of $0$ is spanned by the vector $\mathbf{1} = (1,
\dots, 1).$

\subsection{Complexity of a graph}
An important invariant of an \emph{unweighed} graph $G$ is the number
$\tauG$ of spanning trees  of $G$; it is sometimes called the {\em
  complexity} of $G$.
By Kirckhoff's theorem (\cite{Kirchhoff}),
\begin{equation}
\label{kirckhoff}
n\; \tauG\ =\ \lambda_1\lambda_2\cdot\ldots\cdot\lambda_{n-1}
\end{equation}
and $\tauG$ is equal to
the determinant of any cofactor of the matrix $\Lap$ and
(it is some times called the {\em determinant of Laplacian}).

The definitions for weighed graphs are the essentially the same,
except that
\begin{equation*}
\tauG_f = \sum_{T\in \text{spanning trees of $G$}} \prod_{e\in E(T)}
f(e).
\end{equation*}
This has a natural interpretation in the framework of electrical
circuits, where $f(e)$ is thought of as the \emph{conductance} of
the edge $e$. See \cite{bolobas}.

\subsection{Variational problems}
\label{varspec} The functions we consider are: The bottom nonzero
eigenvalue $\lambda_1$ and
\begin{equation*}
\log \det {\Delta}^* = \sum_{i=1}^{n-1} \log \lambda_i.
\end{equation*}

The bottom non-trivial eigenvalue $\lambda_1$ can be alternatively
defined by the Rayleigh-Ritz quotient:

\begin{equation}
\label{rritz}
\lambda_1 = \min_{\sum_{i=1}^n x_i = 0} \frac{\langle x, \Delta
  x\rangle}{\langle x, x \rangle}
\end{equation}

 From this definition, the concavity of $\lambda_1$ over $P(G)$ is
immediate.

The concavity of $\log\det{\Delta}^*$ is somewhat trickier. First
we show:

\begin{thm}
\label{logconc}
The logarithm of the determinant is a concave function on the set of
positive definite symmetric matrices
\end{thm}

\begin{proof}
Let $Q$ be such a matrix, and
let
$$
Q(t)\ =\ Q+tB,\ \ t\in\reals
$$
be a line of symmetric matrices through $Q$.  Then
$$
\frac{d\log\det (Q(t))}{dt}\ =\ \tr(BQ^{-1}),
$$
and
$$
\frac{d^2\log\det (Q(t))}{dt^2}\ =\ -\tr(BQ^{-1}BQ^{-1}).
$$

It suffices to show that the last trace is strictly positive. The
matrix $R=Q^{-1}$ is positive definite, so can be conjugated by an
orthogonal matrix $P$ to a diagonal matrix $D$, where $D_{ii} >
0$. So, we can rewrite

$$\tr(BRBR) = \tr(BODO^tBODO^t)=\tr((O^tBO)D(O^tBO)D).$$

Let $B'=O^t B O$. $B'$ is still symmetric. We see that $\tr(BRBR)
= \tr(B'DB'D)$. Now, let $d$ be the vector of the diagonal entries
of $D$. It is not hard to check that $\tr(B'DB'D) = d^t
\mathcal{B} d,$ where $\mathcal{B}_{ij} = b_{ij}^2.$  Note,
however, that by our assumptions, all the entries of $d$ are
strictly positive, so  $d^t \mathcal{B} d > 0$, and the result
follows.
\end{proof}

Now we can prove

\begin{thm}
\label{lggconc}
The function $\log \det {\Delta}^*$ is concave on $P(G)$.
\end{thm}

\begin{proof}
The vector $\mathbf{1}=(1, \dots, 1)$ is the zero eigenvector of $\Delta(G_f)$ for
any edge-valuation $f$ in $P(G)$. Thus, the restriction $\Delta ^*$ of
$\Delta$ to the orthogonal complement of the subspace generated by
$\mathbf{1}$ is a symmetric positive-definite operator, whose entries
as a matrix, furthermore, are obviously linear in those of $\Delta$,
no matter which basis of $\mathbf{1}^\perp$ we take. The result now
follows immediately from Theorem \ref{logconc}
\end{proof}

In the sequel, we characterize the extremal valuations for girth,
$\lambda_1,$ and $\log \det {\Delta}^*$ on our deformation spaces $P,
C,$ and $T$.

\section{Maximal girth valuations}
\label{girth2}
\subsection{Maximum in $P(G)$}
Let $G$ be a fixed graph, and suppose that $f_{\max} \in P(G)$ is
such that the $g(f)$ is maximal. There are two, somewhat
different, cases to consider: the first is when $f_{\max}$ is an
interior point of $P(G)$ (\textit{i.e.} no $f(e)$ vanishes), the
second is when $f$ is a boundary point (so that one for one or
more edges $e$, $f(e)= 0$). We will examine the interior point
case first, since it contains the crucial ideas, and is slightly
simpler.

\subsection{Interior maximum} \label{gmax} The idea is that we use something like a
piecewise-linear version of Lagrange multipliers. To wit, suppose
that $f_{\max}$ is our maximal point. That means that there is a
collection of cycles $C_1, \dots, C_n$, such that $\ell(C_i) =
g_f,$ while $\ell(C) > g_f$ for any other cycle $C$. Consider a
small perturbation $g$ of the valuation $f$: $f_1 = f + t h.$
Since $g$ still has to lie in $P(G)$, we must have $\langle h,
\mathbf{1}\rangle = 0.$ The condition that $f$ is maximal is
equivalent to saying that $g_{f_1} \leq g_f$. However, for $t$
sufficiently small, a shortest cycle for the valuation $f_1$ has
to be one of the cycles $C_1, \dots, C_n$, thus the hypothesis
that $g_{f_1} \leq g_f$ means that $\min \ell(C_i)$ at the
weighing $f_1$ has to be smaller than $g_f$. Consider the
quantities $H_i = \sum_{e\in C_i} h(e)$. We know that at least one
of them has to be negative, but this (by multiplying by $-1$ if
necessary) is so if and only if $\exists i, j$, such that $\sgn
H_i = - \sgn H_j$, or else all the $H_i$ vanish. The necessary and
sufficient conditions follow from Farkas' Lemma:

\begin{thm}[Farkas Lemma]\label{fark} 
Let $v_1, \dots, v_n, u \in \mathbb{R}^k$. Then there 
exists a vector $w \in \mathbb{R}^k$, such that $\langle w, v_i 
\rangle \geq 0, \quad 1\leq i \leq n$ (at least one inner product 
being positive) and $\langle u, w\rangle = 0$ if and only if $u$ 
is \emph{not} in the open convex cone generated by the $v_i$. 
\end{thm}

\medskip\noindent
\textbf{Remark.} $u$ is in the convex cone generated by the $v_1,
\dots, v_n$ if there exist $\mu_1, \dots, \mu_n$ either all
negative or all positive, such that $u = \sum_{i=1}^n \mu_i v_i$.

\begin{proof}[Proof of Farkas Lemma]
Suppose first that $$u = \sum_{i=1}^n \mu_i v_i, \qquad \mu_i > 0,
1\leq i \leq n.$$ Take any $w$ such that $\langle w, u\rangle =
0$. Then
\begin{equation*}
0 = \langle w, u \rangle = \sum_{i=1}^n \mu_i \langle w,
v_i\rangle.
\end{equation*}
Since the $\mu_i$ are all positive, not all of the inner products
$\langle w, v_i\rangle$ can be positive, so $w$ does not satisfy
the hypotheses of the theorem.

Suppose now that $u$ is not in the open cone $C$ generated by
$v_1, \dots, v_n$. Consider the projection of $C$ onto the
subspace $u^\perp$ orthogonal to $u$. This is again an open convex
cone $C_u$, which omits at least one point of $u^\perp$ (the
origin). Therefore it is a proper cone, 
and is thus
contained in a half-space $H^+$, and thus the positive normal
vector to $\partial H^+$ has positive inner product with any
vector in the projection of $C$, and hence with any vector in $C$
(since a vector in $C$ can be written as a sum of a vector in
$C_u$ with a multiple of $u$).
\end{proof}

Theorem \ref{fark} can be generalized as follows:

\begin{thm}
\label{fark2}
Let $v_1, \dots, v_n, u_1, \dots, u_m \in\mathbb{R}^k$. 
Then there exists a vector $w \in \mathbb{R}^k$, 
such that $\langle w, v_i \rangle \geq 0, \quad 1\leq i \leq n$ 
(with at one inner product positive) and 
$\langle u_j, w\rangle = 0,\quad 1\leq j\leq m$ 
if and only if \emph{no} linear combination 
$\sum_{j=1}^m a_j u_j$ is in the open convex cone generated by the 
$v_i$. 
\end{thm}

\begin{proof}
If some linear combination $u=\sum_{j=1}^m u_j$ lies in the open
cone $C$, then the same argument as in the beginning of the proof
of Theorem \ref{fark} shows the non-existence of the requisite
$w$. Otherwise, if $u_1, \dots, u_j$ span $\mathbb{R}^k$, there is
nothing left to prove. Assume then that they span a proper
subspace $U$, and project $C$ onto the orthogonal complement, to
get $C_U$. $C_U$ omits the origin by assumption, and the same
argument as in the proof of Theorem \ref{fark} completes the
proof.
\end{proof}

\begin{remark}
\label{strictfark} The above theorems \ref{fark} and \ref{fark2}
do not address the question of when the there exists a nonzero
vector such that the inner products with the $u_j$ \emph{and}
$v_i$ are all zero. This, however, is obviously true if and only
if the span of all of the $v_i$ together with all of the $u_j$ is
a proper subspace of $\mathbb{R}^k$.
\end{remark}

Theorems \ref{fark} and \ref{fark2} and Remark \ref{strictfark}
combine to give the following characterization of the extremal
points of girth in $P(G)$, $T(G)$ and $C(G)$, which we state in
the Theorem \ref{omnibus} below. First

\medskip\noindent
\textbf{Notation.} The systoles of $G$ corresponding to a weighing
$f$ are cycles $s_1, \dots, s_k$ whose length is equal to the
girth of $G$ with the weighing $f$. We  call \emph{edge systoles}
the vectors $\mathfrak{s}_1, \dots, \mathfrak{s}_k$ in
$\mathbf{R}^{E(G)}$ whose $e$-th coordinate is $1$ if $e$ is
contained in the corresponding cycle $s_j$. We call the
\emph{vertex systole} corresponding to $s_i$, the vector
$\sigma_i$ in $\mathbb{R}^{V(G)}$, whose $v$-th coordinate is $1$
if $v$ is incident to $s_i$, and $0$ otherwise. The \emph{vertex
vector} $w_v$ is the vector in $\mathbb{R}^{E(G)}$ whose $e$-th
coordinate is $0$ unless $e$ is incident to the vertex $v$, in
which case the coordinate is $1$. The \emph{degree vector} $d(G)$
is the vector in $\mathbb{R}^{V(G)}$ whose $v$-th coordinate is
the degree of the vertex $v$.
\begin{thm}
\label{omnibus} A weighting $f\in P(G)$ is maximal for girth if
and only if the constant vector $\mathbf{1}$ lies in the open cone
generated by the edge systoles of $G$ with the weighting $f$. The
maximal weighing $f$ is unique if and only if the edge systoles of
$G$ corresponding to the weighing $f$ together with the constant
vector span the whole space $\mathbb{R}^{E(G)}$.

A weighing $f$ in $T(G)$ is maximal for girth if and only if some
linear combination of the vertex vectors $w_1, \dots, w_{V_G}$
lies in the open cone generated by the edge systoles of $G$ with
the weighing $f$. The maximal weighing $f$ is unique if and only
if the edge systoles and the vertex vectors span
$\mathbb{R}^{E(G)}$.

A weighing $f$ in $C(G)$ is maximal for girth if and only if the
degree vector $d(G)$ is contained in the open cone generated by
the vertex systoles of $G$. The maximal weighing is unique if and
only if the degree vector together with the vertex systoles span
$\mathbb{R}^{V(G)}.$
\end{thm}

\section{The tree number}
\label{maxtreesec} By the weighted version of Kirckhoff's theorem
(\cite{bolobas})
\begin{equation}\label{det:tree-exp}
    \tauG\ =\ \sum_{T\in\cT(G)} \; \prod_{e_j\in T} x_j
\end{equation}
where the sum is taken over the set $\cT(G)$ of the spanning trees
of $G$.

We will find necessary and sufficient condition for a valuation
$f$ to be the critical point for $\tauG$ (which is the same as
being maximal by $\log \det \Lap^*$, by the discussion in the
Introduction) on $P(G)$, $C(G)$ and $T(G)$. It should be noted
that such a critical point might not exist, and we might have to
look for boundary maxima. Our methods can be easily adapted to
deal with those cases as well, and since writing down the
conditions is somewhat more cumbersome, we leave this to the
reader.   

\subsection{Maximum in $P(G)$.}
We start with $P(G)$, since the result in that case is the
simplest to state, and seems, at least at the moment to have the
simplest combinatorial interpretation.
  Finding the maximum of $\tauG$ on $P$
is a Lagrange multiplier problem. The condition for $x\in P(G)$ to
be a critical point for $\tauG$ is
\begin{equation}
\frac{\partial\tauG}{\partial x_1}\ =\
\frac{\partial\tauG}{\partial x_2}\ =\ \ldots\ =\
\frac{\partial\tauG}{\partial x_m}
\end{equation}

The partial derivatives above are given by $$ \tau_j\ =\
\frac{\partial\tauG}{\partial x_j}\ =\ \sum_{e_j\in T\in\cT(G)}
\prod_{k\neq j;e_k\in T} x_k. $$ The ratio $\tau_j/\tauG$ is
called the {\em effective resistance} of $e_j$.

We have thus proved:
\begin{prop}\label{prop:detcrit}
The graph valuation $f$ is maximal for $\tauG$ in $P(G)$  if and
only if the effective resistances of all edges are the same.
\end{prop}

If an unweighted graph satisfies the assumptions of Proposition
\ref{prop:detcrit} then every edge of this graph is contained in
{\em the same} number of spanning trees.  Such graphs were studied
by Godsil in \cite{Godsil}; he calls these graphs {\em
equiarboreal}.  Obviously, all edge-transitive graphs (the
automorphism group acts transitively on the edges) are
equiarboreal.\footnote{See \cite{Bouwer} for examples of
edge-transitive graphs which are not vertex-transitive.} Godsil
gives several more sufficient conditions for a graph to be
equiarboreal; in particular, any distance-regular graph and any
color class in an association scheme is equiarboreal (the least
restrictive condition Godsil gives is for a graph to be {\em
1-homogeneous}).  By an easy counting argument one can show that
for an unweighted equiarboreal graph
\begin{equation}\label{equiarb:effres}
   T_1\ =\ T_2\ =\ \ldots\ =\ \tauG\cdot (n-1)/m, 
\end{equation}
where $T_j$ is the number of spanning trees containing $e_j$ 
(this is actually the result of Foster, cf. \cite{Foster}) so the
necessary condition for a graph to be equiarboreal is that $m$
divide $(n-1)\tauG$.
\begin{remark}
Any tree is equiarboreal.
\end{remark}

We remark that the graphs which have the most spanning trees among
the regular graphs with the same number of vertices are not
necessarily equiarboreal, and vice versa.  For example, the
$8$-vertex M{\"o}bius wheel (cf. \cite{Biggs:book}) which has the
most spanning trees among the $8$-vertex cubic graphs is not
equiarboreal (cf. also \cite{Valdes}), while the {\em cube} (which
is certainly edge-transitive, hence equiarboreal) has the {\em
second biggest} number of spanning trees among the $8$-vertex
cubic graphs.

\subsection{Maxima in $T(G)$ and $C(G)$}
The Lagrange multiplier method of the previous section works just
as well in $T(G)$ and $C(G)$. We leave the (easy) computation to
the reader, and just summarize the results in

\begin{thm}
\label{othermax} A valuation $f$ is maximal in $T(G)$ if and only
if there exists constants $\lambda_1, \dots, \lambda_{V(G)}$, such
that if the edge $e$ has endpoints $v_i$ and $v_j$, then
\begin{equation*}
\tau(e) = \lambda_i + \lambda_j.
\end{equation*}

A valuation $f$ is maximal in $C(G)$ if and only if for any two
vertices $v$ and $w$
\begin{equation*}
\deg w \sum_{\text{$e$ incident to $v$}} \tau(e) = \deg v
\sum_{\text{$e$ incident to $w$}} \tau(e).
\end{equation*}
\end{thm}

If we ask the same question as previously -- when is the constant
valuation maximal? -- the condition for a maximum in $T(G)$ does
not appear to have an obvious combinatorial interpretation. The
condition for the maximum in $C(G)$ can be restated in the
following way:

\begin{cor}
Let $d_T(v) = \sum_{\text{spanning trees $T$}} \text{$\deg v$ in
$T$}$. Then, if the constant valuation is maximal for $\tauG$ on
$C(G)$, then for any two vertices $v$ and $w$,
\begin{equation*}
\frac{d_T(v)}{\deg v} = \frac{d_T(w)}{\deg w}.
\end{equation*}
\end{cor}

\section{Eigenvalues of the Laplacian}
\label{maxevsec} To find the condition for maximality with respect
to the bottom non-zero eigenvalue of the Laplacian, we will use
the Rayleigh-Ritz characterization of of $\lambda_1$. This implies
immediately that:

\begin{thm}
\label{lambdacond}
Let $f$ be the weighing on $G$
(in our application, $S$ could be
any one of $P(G)$, $T(G)$, $C(G)$, but it could be anything). Let
$E_{\lambda_1}$ be the eigenspace corresponding to $\lambda_1$. Let
$Q$ be any infinitesimal variation (that is, an element of the tangent
space of $S$) of the valuation, and $Q_\Delta$
the induced variation of the Laplacian matrix. Then the quadratic form
given by $Q_\Delta$ restricted to $E_{\lambda_1}$ is indefinite if and
only if $f$ is maximal with respect to $\lambda_1$.
\end{thm}

\begin{proof}
The argument is a version of that given in the beginning of section
\ref{gmax}. We use the Rayleigh-Ritz quotient characterization (given
in eq. \ref{rritz}). The space $E_{\lambda_1}$ is precisely the set of
vectors where the minimum is attained, so at any unit vector $x \notin
E_{\lambda_1}$, $\langle x, (\Delta + t Q_\Delta) x \rangle$ is
strictly greater than $\langle y, (\Delta + t Q_\Delta) y \rangle$ for
$y$ a unit vector in $E_{\lambda_1}$, for $t$ sufficiently
small. Thus, the first variation of $\lambda_1$ is given by the first
variation of $\lambda_1$ restricted to $E_{\lambda_1}$, and that is
given precisely by the restriction of the quadratic form given by
$Q_\Delta$. Now, if that were definite, we would be
able to increase $\lambda_1$ by applying either the variation
$Q_\Delta$ or $-Q_\Delta.$
\end{proof}

Note now that the space of all possible variations of the Laplacian
induced by changes in the edge valuations has a natural linear
structure (one can think of it as a subspace of the tangent space to
symmetric matrices). Call that space $V_{\text{var}}.$ If $M \in
V_{\text{var}}$, then $x^t M x$ can be thought of as a scalar product
of $M$ with a vector $P_x$, whose $ij$-th coordinate is given by $x_i
x_j$ (this is just the outer product of $x$ with itself, the letter
$P$ is used to point out that when $x$ is a unit vector, $P_x$ is just
the projection on the subspace generated by $x$). Let
$\mathcal{P}_{\lambda_1} = \left\{P_x \quad | \quad x \in
  E_{\lambda_1}\right\}$. If
$S^\perp$ is the orthogonal complement to the tangent space of the
deformation $S$, Theorem \ref{fark2} (whose proof does not use the
finiteness of the sets involved) gives us:

\begin{thm}
\label{fundlam}
A valuation $f$ is maximal in $S$ with respect to $\lambda_1$ if and
only if the intersection of $S^\perp$ with the open cone generated by
$\mathcal{P}_{\lambda_1}$ is nonempty.
\end{thm}

What is the ``open cone generated by $\mathcal{P}_{\lambda_1}$''~? It
is an easy exercise to show that this is precisely the set of
positive self-adjoint operators on $E_{\lambda_1}$ (that is, operators
for which $E_{\lambda_1}$ is an invariant subspace~; which are
positive on that subspace, and zero elsewhere).
 so Theorem \ref{fundlam} can be restated as:
\begin{thm}
\label{fundlam1}
A valuation $f$ is maximal in $S$ with respect to $\lambda_1$ if and
only if $S^\perp$ contains a positive self-adjoint operator $\Lambda$
on $E_{\lambda_1}$.
\end{thm}

\begin{cor}
\label{fundcor}
If $E_{\lambda_1}$ is one-dimensional, then $f$ is maximal if and only
if $S^\perp$ is spanned by $P_v$, where $v$ is a unit eigenvector of
$\lambda_1$.
\end{cor}

All the above might sound somewhat abstract, so let us now
specialize to the the deformation spaces we have in mind. First,
consider $P(G)$. In this case, it is easy to check that
\begin{equation*}
\frac{\partial \Delta}{\partial f(e)} = Q_e,
\end{equation*}
where, if the endpoints of $e$ are $v_i$ and $v_j$, then $Q_{ii} =
Q_{jj} = 1$~; $Q_{ij}=Q_{ji} = -1$, and all of the other entries are
$0$.
The general variation of $\Delta$ is given by $Q_{\mathbf{\alpha}} =
\sum_e \alpha_e Q_e,$ and in order to stay in $P(G)$, we must have
$\langle \mathbf{\alpha}, \mathbf{1} \rangle = 0.$

It can be seen that the variation space $S$ of the Laplacians is spanned
by the vectors $Q_{e_0} - Q_e$, where $e_0$ is an arbitrary fixed
edge.

\subsection{The first eigenvalue $\lambda_1$ appears without
  multiplicity.}
If the eigenspace of $\lambda_1$ is one-dimensional, and the
eigenvector is $v$, then by Corollary \ref{fundcor}, $v^t (Q_{e_0}
- Q_e) v = 0$, for all $e$. If $e$ is an edge with endpoints $x$ and
$y$, then a calculation shows that $v^t Q_e v = (v_x - v_y)^2$, and so for
a maximal valuation, we must have
\begin{equation}
\label{df:const}
v_x - v_y = \pm c
\end{equation}
for \textit{any} adjacent pair of vertices $x, y$. We study graphs
which have an eigenvector satisfying the condition given by eq.
(\ref{df:const}) in section \ref{constgrad}, but it is \textit{a
priori} clear that this condition is very rarely satisfied, and
``usually'' graphs maximal for $\lambda_1$ have a
higher-dimensional first eigenspace. Curiously, the same holds for
the (smaller) deformation spaces $T(G)$ and $C(G)$. Indeed,
consider first $T(G)$. There, the deformation space of the
Laplacians is spanned by matrices $Q_{e_1} - Q_{e_2}$, where $e_1$
and $e_2$ have a vertex in common. Thus, the same computation as
that leading to eq. (\ref{df:const}) gives that the eigenvector of
$\lambda_1$ for a critical graph must satisfy:

\begin{equation}
\label{dfconst2}
v_x - v_y = \pm c_x,
\end{equation}
for any adjacent pair of vertices $x, y$. \textit{A priori,} this
seems somewhat weaker than the condition (\ref{df:const}) (since $c_x$
now depends on $x$), but in fact it is clear that for a
\textit{connected} graph $G$, it is equivalent~; the case of $G$
disconnected is different, but not particularly interesting.

For $C(G)$, the deformation space of Laplacians is generated by
the differences $\deg w M_v - \deg v M_w$, where $M_i$ is the
matrix whose $ii$-th entry is the degree of the $i$-th vertex;
$M_{ij}$ is equal to $-1$ if $v_j$ is incident to $v_i$, likewise
$M_{ji}$, and all other $M_{jk}$ are equal to $0$. If $v$ is a
vector, then $$v^t M_x v = \deg x v_x^2 - 2 \sum_{y \sim x} v_x
v_y = 2 v_x \sum_{y \sim x}(v_x - v_y) - \deg x v_x^2.$$ If $v$ is
an \textit{eigenvector} of $G$ with eigenvalue $\lambda$, then
$\sum_{y \sim x} (v_x - v_y) = \lambda v_x,$ and so $$v^t M_x v =
(2 \lambda - \deg x)v_x^2.$$ From the equation $v^t (\deg w M_v -
\deg v M_w) v = 0$, it follows that:
\begin{equation}
\label{confeq}
(\frac{2 \lambda}{\deg x} - 1) v_x^2 = (\frac{2 \lambda}{\deg y} - 1) v_y^2.
\end{equation}
In particular, note that when the graph $G$ is regular, it follows
that
\begin{equation}
\label{confreg}
|v_x| = |v_y|,
\end{equation}
for any two vertices $x$, $y$.

\medskip\noindent
{\bf Remark.} It is not difficult to construct regular graphs
which have an eigenvector satisfying eq. (\ref{confreg}): any such
graph is constructed by taking an $l$ regular bipartite graph,
whose vertex set is the union of the sets $R$ of red vertices and
$B$ of black vertices, and constructing $k$-regular graphs with
vertex sets $R$ and $B$ respectively (then adjoining their edge
sets to that of the original bipartite graph). Then the function
which is $1$ on $R$ and $-1$ on $B$ is an eigenvector with
eigenvalue $2 l$. It is much less clear that this can be done in
such a way that $2 l$ is \textit{the
  lowest} eigenvalue.

\subsection{The general case.} When the eigenspace of $\lambda_1$ has
dimension possibly greater than $1$, Theorem \ref{fundlam1},
compact with the finite-dimensional spectral theorem (that a
positive self-adjoint operator can be diagonalized, with respect
to an orthonormal basis, with positive weights) gives us the
following extensions of the results of the previous subsection:

\begin{thm}
\label{multidim1}
In order for a valuation $f$ to be maximal for $\lambda_1$ with
respect to $P(G)$, it is necessary and sufficient for there to be an
orthogonal basis $v_1, \dots, v_d$ of $E_{\lambda_1}$, and a
collection of non-negative constants $c_1, \dots, c_d$,  not all zero,
and a constant $c > 0$ such that for any pair of adjacent vertices $x,
y$ of $G$
\begin{equation}
\label{multeq1}
\sum_{i=1}^d c_i (v_i(x) - v_i(y))^2 = c.
\end{equation}

In order for $f$ to be maximal for $\lambda_1$ with respect to $T(G)$,
the same condition (\ref{multeq1}) holds, assuming that $G$ is
connected.

In order for $f$ to be maximal for $\lambda_1$ with respect to $C(G)$,
there must be constants as above, such that for any vertex $x$ of $G$,
\begin{equation}
(\frac{2 \lambda}{\deg x} - 1)\label{multeq2} \sum_{i=1}^d (c_i
v_{i}^2(x)) = c.
\end{equation}
\end{thm}

\begin{remark} The condition \eqref{multeq1} gives an embedding 
of the edge set $E(G)$ into a $(d-1)$-dimensional \textit{ellipsoid} 
$$
\sum_{i=1}^dc_iz_i^2=c
$$ 
by the ``differentials'' $z_i=dv_i(e)=v_i(x)-v_i(y)$, where we have 
chosen an arbitrary orientation of the edge $e=(xy)$.  The corresponding 
vertex condition gives, for a regular graph, a similar embedding 
of the vertex set $V(G)$ by the eigenvectors $v_i$.  
\end{remark}

\section{Graphs with an eigenvector of constant gradient}
\label{constgrad}
We now study connected graphs which admit an eigenvector $f:V\to\reals$
satisfying \eqref{df:const} for some $c\geq 0$. If $c=0$ then $f$
is a multiple of a constant vector and so has eigenvalue zero
which is a contradiction. If $c\neq 0$ then it is easy to see that
the graph $G$ cannot have odd cycles and hence is bipartite.
Namely, let $u_1u_2\ldots u_l$ be a cycle.  Then (putting
$u_l=u_0$) $\sum_{i=1}^l (f(u_i)-f(u_{i-1}))=0$.  But each term in
the sum is equal to $\pm c$, and since the number of terms in the
sum is odd, they cannot add up to $0.$

We now want to study the unweighted $k$-regular graphs which have
an eigenvector (corresponding to an eigenvalue $\mu>0$) satisfying
\eqref{df:const} (without necessarily assuming that $\mu$ is
simple). We shall rescale the eigenvector so that $c=1$ in
\eqref{df:const}.  From  \eqref{df:const} it follows that for each
vertex $u$ the expression $\mu\cdot f(u)$ can only take one of the
values $k, k-2,k-4,\ldots, -k+2,-k$. Consider first the vertex
$u_0$ where $f(u)$ takes its maximal value $a$ (by changing the
sign if necessary we can assume that $a>0$). It follows that $f$
takes value $a-1$ on all the neighbors of $u_0$, hence $$
   a\mu\ =\ k
$$

Next, consider any neighbor $u_1$ of $u$.  The value of $f$ at any
neighbor of $u_1$ can be either $a$ (let there be $r_1\geq 1$ such
neighbors; $u_0$ is one of them); or $a-2$ (it follows that there
are $k-r_1$ such neighbors).  From the definition of the Laplacian
it follows that $$
   \mu(a-1)\ =\ k-2r_1
$$ It follows from the last two formulas that
\begin{equation}\label{eig:even}
   \mu\ =\ 2r_1
\end{equation}
where $r_1\geq 1$ is a positive integer.  If $r_1=k$, then
$\mu=2k$ is the largest eigenvalue of $\Lap$.

We next define the {\em level} of a vertex $u$ to be equal to $j$
if $f(u)=a-j$; we denote the set of all vertices of $G$ at level
$j$ by $G_j$.  It is easy to see that if $u\in G_j$ has $r_j$
neighbors where $f$ takes value $a-j+1$ then $$
   \mu(a-j)\ =\ k-2r_j
$$ It follows that $r_j$ is the same for all $u\in G_j$. Using
\eqref{eig:even} we see that $$
   r_1\cdot j\ =\ r_j
$$ Consider now a ``local minimum'' $u\in G_N$.   Then $r_N=k$,
and we see that
\begin{equation}\label{eig:divides}
   r_1\; |\; k
\end{equation}

Let $n_j$ denote the number of vertices in $G_j$.  Counting the
vertices connecting $G_j$and $G_{j+1}$ in two different ways, we
see that for all $0\leq j\leq N-1$, $$ n_j(k-r_j)\ =\
n_{j+1}r_{j+1} $$ Consider the case $r_1=1,\mu=2$.  It follows
from the previous calculations that $r_j=j$ and that $N=k$.
Accordingly, $n_j=n_0{k\choose j}$ and
\begin{equation}\label{vertn:dfconst:2}
|G|\ =\ 2^k\; n_0
\end{equation}
We next describe a class of graphs admitting an eigenvector of
$\Lap$ with $\mu=2$ satisfying \eqref{df:const}.

An obvious example of such a $k$-regular graph is the $k$-cube,
and any such graph has the same number of vertices as a disjoint
union of $n_0$ cubes by \eqref{vertn:dfconst:2}.  Start now with
such a union, choose the partition of the vertices of each cube
into ``levels'' and take two edges $u_1u_2$ and $u_3u_4$ in two
different cubes such that $u_1,u_3$ are both in level $j$ while
$u_2,u_4$ are both in level $j+1$. If we perform an edge switch $$
(u_1u_2),(u_3u_4)\ \to\ (u_1u_4),(u_3u_2) $$ then the number of
the connected components of our graph will decrease while the
eigenvector $f$ will remain an eigenvector with the same
eigenvalue.

Performing sequences of edge switches as described above, we
obtain examples of connected graphs satisfying \eqref{df:const}
and \eqref{vertn:dfconst:2} for any $n_0$.  Conversely, it is easy
to show that starting from a graph satisfying \eqref{df:const} and
\eqref{vertn:dfconst:2} and having chosen a partition of its
vertices into levels one can obtain $n_0$ disjoint $k$-cubes by
performing a sequence of edge switches as above.

We now want to consider the case when $\mu=\mu_1$ is the lowest
eigenvalue of the Laplacian.  The first remark is that then
necessarily $\mu\leq k$, and $\mu=k$ only if $G=K_{k,k}$.  Next,
we want to consider ``small'' $k$ for which $k-2\sqrt{k-1}$ (the
``Ramanujan bound'') is less than $2$ (this happens for $3\leq
k\leq 6$). It then follows from the results of Alon (\cite{Nilli})
that the diameter of $G$ (and hence the number of vertices in $G$)
is bounded above.
\begin{prop}\label{mu1:finite}
For $3\leq k\leq 6$ there are finitely many $k$-regular graphs for
which the condition \eqref{df:const} is satisfied for an
eigenvector of $\mu_1$.
\end{prop}

We next discuss graphs which have an eigenvector satisfying
\eqref{df:const} with the eigenvalue $\mu=2r_1>2$.  Recall that by
\eqref{eig:divides} $r_1|k$. By counting the edges connecting the
vertices in two consecutive levels one can show (as for $\mu=2$)
that the number of vertices satisfies
\begin{equation}\label{vertn:general}
|G|\ =\ 2^{(k/r_1)}\; n_0
\end{equation}
Also, since any vertex $u_1\in G_1$ has $r_1$ distinct neighbors
in $G_0$, $$ n_0\geq r_1. $$ It is easy to construct examples of
regular graphs which have eigenvectors with the eigenvalue $\mu>2$
satisfying \eqref{df:const}; the construction is similar to that
for $\mu=2$.

We summarize the previous results:
\begin{thm}\label{eigencrit:character}
Let $G$ be a $k$-regular graph which has an eigenvector of $\Lap$
with an eigenvalue $\mu$ satisfying \eqref{df:const}.  Then $G$ is
bipartite, $\mu=2l$ is an even integer dividing $2k$, the number
of vertices of $G$ is divisible by $2^{(k/l)}$, and for $n_0\geq
l$ there exist such graphs with $n=2^{(k/l)}n_0$ vertices.
\end{thm}


\begin{thebibliography}{}

\bibitem[Big93]{Biggs:book}
N. Biggs.
\newblock {\em Algebraic graph theory (2nd ed).}
\newblock Cambridge Univ. Press, 1993.

\bibitem[Bol98]{bolobas}
B. Bollob{\'a}s. \textit{Modern Graph Theory\/,} Springer Verlag,
New York, 1998.

\bibitem[Bou]{Bouwer}
I. Bouwer.
\newblock {\em On edge but not vertex transitive regular graphs.}
\newblock J. Comb. Th. B, 12:32--40, 1972.

\bibitem[Fie]{Fiedler}
M. Fiedler.
\newblock {\em Some minimax problems for graphs.}
\newblock Discr. Math, 121:65--74, 1993.

\bibitem[Fos]{Foster}
R. Foster.
\newblock {\em The average impedance of an electrical network.}
\newblock Contib. to Applied. Mechanics (Reissner Ann. Volume),
Edwards Bros, 333--340, 1949.

\bibitem[God81]{Godsil}
C. Godsil.
\newblock {\em Equiarboreal graphs.}
\newblock Combinatorica, 1:163--167, 1981.

\bibitem[Kir]{Kirchhoff}
F. Kirchhoff.
\newblock {\em {\"U}ber die Aufl{\"o}sung der Gleichungen, auf welche man
bei der Untersuchung der linearen Verteilung galvanischer Str{\"o}me
gef{\"u}rt wird.}
\newblock Ann. Phys. Chem. 72:497--508, 1847.

\bibitem[Nil]{Nilli}
A. Nilli.
\newblock {\em On the second eigenvalue of a graph.}
\newblock Discr. Math, 91:207--210, 1991.

\bibitem[OPS]{OPS}
B. Osgood, R. Phillips and P. Sarnak.
\newblock  {\em Extremals of Determinants of Laplacians.}
\newblock J. Func. Anal, 80:148--211, 1988.


\bibitem[Riv99]{Riv99}
I.~Rivin.
\newblock {\em Growth in free groups (and other stories)}.
\newblock xxx.lanl.gov preprint math.CO/9911076.

\bibitem[Sch]{Schrijver}
A. Schrijver.
\newblock {\em Theory of linear and integer programming.}
\newblock John Wiley \& Sons, 1990.

\bibitem[Val]{Valdes}
L. Valdes.
\newblock {\em Extremal properties of spanning trees in cubic graphs.}
\newblock Congr. Numer, 85:143--160, 1991.

\end{thebibliography}
\end{document}